\documentclass[a4paper,10pt]{amsart}
\usepackage[utf8]{inputenc}
\usepackage{amsthm}
\usepackage{amsmath}
\usepackage{bm}
\usepackage{amsfonts}
\usepackage{amssymb}
\usepackage{mathtools}
\usepackage{mathrsfs}
\usepackage{setspace}
\usepackage{xcolor}
\usepackage{textgreek}
\usepackage{todonotes}
\usepackage{thmtools} % solves problem of shared counters in statements and autoref

\usepackage[pdfdisplaydoctitle,colorlinks,breaklinks,urlcolor=blue,linkcolor=blue,citecolor=blue]{hyperref} %must always be the last

\newcommand{\A}{\mathcal{A}}

\newcommand{\cyl}{\mathcal{C}}
\newcommand{\D}{\mathcal{D}}

\newcommand{\E}{\mathcal{E}}

\newcommand{\F}{\mathcal{F}}
\newcommand{\G}{\mathcal{G}}

\renewcommand{\L}{\mathcal{L}}

\newcommand{\N}{\mathbb{N}}

\newcommand{\PP}{\mathbb{P}}

\newcommand{\R}{\mathbb{R}}
\newcommand{\RR}{\mathcal{R}}
\renewcommand{\S}{\mathcal{S}}

\newcommand{\Z}{\mathbb{Z}}

\DeclareMathOperator{\id}{Id}
\let\div\relax
\DeclareMathOperator{\div}{div}

\renewcommand{\epsilon}{\varepsilon}

\newcommand{\one}{\bm{1}}

\DeclareMathOperator{\curl}{curl}
\let\div\relax
\DeclareMathOperator{\div}{div}

\newcommand{\de}{\partial}

\newcommand{\set}[1]{\left\{#1\right\}}
\newcommand{\pa}[1]{\left(#1\right)}
\newcommand{\bra}[1]{\left[#1\right]}
\newcommand{\abs}[1]{\left|#1\right|}
\newcommand{\norm}[1]{\left\|#1\right\|}
\newcommand{\brak}[1]{\left\langle#1\right\rangle}

\newcommand{\expt}[2][]{\mathbb{E}_{#1}\left[#2\right]}
\def\XXint#1#2#3{{\setbox0=\hbox{$#1{#2#3}{\int}$}
		\vcenter{\hbox{$#2#3$}}\kern-.5\wd0}}
\newtheorem{thm}{Theorem}[section]

\newtheorem{definition}[thm]{Definition}

\newtheorem{lemma}[thm]{Lemma}

\theoremstyle{remark}
\newtheorem{rmk}[thm]{Remark}

\numberwithin{equation}{section}

\title[Invariant Measures 2D PE]{Gaussian Invariant Measures\\ and Stationary Solutions\\of 2D Primitive Equations}
\author[F. Grotto]{Francesco Grotto}
  \address{Scuola Normale Superiore, Piazza dei Cavalieri, 7, 56126 Pisa, Italia}
  \email{\href{mailto:francesco.grotto@sns.it}{francesco.grotto@sns.it}}
\author[U. Pappalettera]{Umberto Pappalettera}
  \address{Scuola Normale Superiore, Piazza dei Cavalieri, 7, 56126 Pisa, Italia}
  \email{\href{mailto:umberto.pappalettera@sns.it}{umberto.pappalettera@sns.it}}
\keywords{invariant measure, primitive equations, regularization by noise, equilibrium statistical mechanics}
\date\today

\begin{document}

\begin{abstract}
We introduce a Gaussian measure formally preserved by the 2-dimensional 
Primitive Equations driven by additive Gaussian noise.
Under such measure the stochastic equations under consideration are singular:
we propose a solution theory based on the techniques developed by Gubinelli and Jara in \cite{GuJa13}
for a hyperviscous version of the equations.
\end{abstract}

\maketitle

%%%%%%%%%%%%%%%%%%%%%%%%%%%%%%%%%%%%%%%%%%%%%%%%%%%%%%%%%%%%%%%%%%%%%%%%%%%%%%%%%%%%%%%%%%%%%
\section{Introduction}\label{sec:introduction}

Primitive Equations constitute a fundamental model in geophysical fluid dynamics.
The present work is devoted to the study of Gaussian invariant measures in the stochastically forced 2-dimensional case:
the model under analysis is thus a stochastic PDE of the form:
\begin{equation}\label{eq:generalmodel}
	\begin{cases}
	\de_t v + v \de_x v + w \de_z v + \de_x p = \D(\Delta)v + \eta, \\
	\de_z p = 0, \\
	\de_x v + \de_z w = 0,
	\end{cases}
\end{equation}
where $(x,z)$ are coordinates of the bounded domain $D=[0,2\pi]^2$
on which suitable boundary conditions are imposed,
$(v,w)$ are the components of the \emph{velocity} vector field,
$p$ is the \emph{pressure},
the term $\D(\Delta)$ describes a dissipation mechanism 
and $\eta$ is a Gaussian stochastic process.

It is in fact the case with $\D(\Delta)=\nu\Delta$, and $\eta=0$ to be usually 
referred to as \emph{2-dimensional Primitive Equations} (2dPE), together with its variants including effects 
such as density and temperature variations, and other geophysical effects. 
When those physical phenomena are neglected,
equations \eqref{eq:generalmodel} have many aspects in common with the 2-dimensional Navier-Stokes equations:
\begin{equation*}
\begin{cases}
\de_t u + (u\cdot \nabla)u+\nabla p=\nu \Delta u,\\
\div u=0,
\end{cases}
\end{equation*}
describing the evolution of the velocity field $u$ of an incompressible viscous fluid.
This familiarity naturally leads to look for applications of concepts and techniques developed
in the extensive theory of Navier-Stokes equations, especially in the 2-dimensional setting.
However, the nonlinearity of 2-dimensional Primitive Equations is in fact harder to treat.

When considering the stochastically forced case, it is well-known that 2-dimensional 
stochastic Navier-Stokes equations (SNS) in their vorticity form:
\begin{equation*}
	\begin{cases}
	\de_t \omega + (u\cdot \nabla)\omega=\nu \Delta \omega+\sqrt{2\nu} \xi,\\
	\div u=0, \, \omega=\curl u,
	\end{cases}
\end{equation*}
$\omega$ being the scalar \emph{vorticity} field, preserve the so-called \emph{enstrophy measure}
when driven by space-time white noise $\xi$.
The enstrophy measure is a Gaussian random distribution, corresponding to space white noise
at the level of $\omega$: its name is due to the fact that its covariance is the quadratic form associated
to the quadratic observable $\norm{\omega}^2_{L^2}$, known as \emph{enstrophy}.
Indeed, enstrophy is a first integral of motion in the case $\nu=0$, the 2-dimensional Euler equations,
and enstrophy measure is the unique, ergodic invariant measure of the linear part of the dynamics
when $\nu>0$.
Notwithstanding the low space regularity under enstrophy measure,
existence and pathwise uniqueness of stationary solutions for SNS in this setting 
are by now classical results due to \cite{DPDe02,AlFe04}.

There exists of course another quadratic invariant for Euler equations:
the energy $\norm{u}^2_{L^2}$. When SNS is driven by space-time white noise at the level of velocity,
the \emph{energy measure}, a white noise at the level of $u$, is \emph{formally} preserved.
The cursive is here in order, because the energy measure regime is so singular that
no solution theory is yet available in this case. Nonetheless, the existing stochastic analysis techniques 
allow to deal with such regime in \emph{hyperviscous} cases, that is replacing the viscous term
$\Delta u$ with $-(-\Delta)^\theta u$, $\theta>1$. Indeed, 
a procedure known as \emph{It\=o trick} in the literature related to regularisation by noise
is employed in \cite{GuJa13} to give meaning and solve SNS under energy measure
with sufficiently strong hyperviscosity. We also mention the recent development \cite{GuTu19},
in which Kolmogorov equations are solved by means of Gaussian analysis tools, 
broadening the result to solutions absolutely continuous with respect to energy measure.

The analogue of vorticity field for 2-dimensional Primitive Equations is $\omega=\de_z v$, 
as the quadratic observable $\norm{\de_z v}^2_{L^2}$ is a first integral
of the 2-dimensional \emph{hydrostatic Euler equations}:
\begin{equation}\label{eq:2dhe}
	\begin{cases}
	\de_t v + v \de_x v + w \de_z v + \de_x p =0, \\
	\de_z p = 0, \\
	\de_x v + \de_z w = 0.
	\end{cases}
\end{equation}
Prescribing the correct additive Gaussian noise $\eta$, the linear part of \eqref{eq:generalmodel}
with $\D(\Delta)=\Delta$ preserves the Gibbsian measure associated to $\norm{\de_z v}^2_{L^2}=\norm{\omega}^2_{L^2}$,
formally defined by
\begin{equation*}
	d\mu(\omega)=\frac1Z e^{-\frac12 \norm{\omega}^2_{L^2}} d\omega,
\end{equation*}
that is, white noise distribution for $\omega$.
However, as we will detail below,
the stationary regime with $\mu$-distributed marginals for \eqref{eq:generalmodel} is
not comparable to the enstrophy measure stationary regime of SNS, because of the
more singular nonlinearity. Indeed, unlike in \cite{DPDe02,AlFe04},
the nonlinear terms of \eqref{eq:generalmodel} can not be defined as distributions when $\de_z v$ has law $\mu$.
Still, as in the case of energy measure SNS, hyperviscosity allows to apply the techniques of \cite{GuJa13}.

In this work we present a solution theory of 2-dimensional Primitive Equations in the hyperviscous setting
$\D(\Delta)=-(-\Delta)^\theta$, for large enough $\theta$ and a suitable stochastic forcing.
The regularising effect of hyperviscosity for Navier-Stokes and Primitive Equations
is well-understood in the deterministic setting, and it is often used in numerical simulations \cite{LJTN11};
we refer to \cite{La58,Li59,LiTeSh92a} and, more recently, \cite{Hu20} for a thorough discussion.
The main contribution of the present paper is thus to introduce a Gaussian invariant measure
in the context of 2-dimensional Primitive Equations, and then to exploit the techniques
of \cite{GuJa13} to provide a first well-posedness result for this singular SPDE in a hyperviscous setting.

Although stochastic versions of Primitive Equations both in two and three dimensions
have already been considered, to the best of our knowledge the existing literature is limited
to regimes more regular than ours. To mention a few relevant previous works,
in \cite{GHZi08,GHTe11,TM10}, 2-dimensional Primitive Equations
are considered with a multiplicative noise taking values in function spaces,
the same is done in the 3-dimensional case in \cite{DeGHTeZi12,GaSu16}, 
and in \cite{GHKuViZi14} the authors prove the existence of an invariant measure in this setting.
In the 2-dimensional cases, large deviation principles are studied in \cite{GaSu12,SuGaLi18}. 
Let us also mention the works \cite{PeTeWi04,BrKaLe05} on deterministic 2-dimensional
Primitive Equations, whose study began with \cite{LiTeSh92a,LiTeSh92b,LiTeSh93},
and \cite{MaWo12} on their inviscid version, by which the vorticity formulation
we present below for our model is inspired.

The paper is structured as follows: in \autoref{sec:vortform} we rigorously introduce a
stochastic version of 2-dimensional Primitive Equations in terms of vorticity $\omega=\de_z v$
and a Gaussian measure formally preserved by the dynamics. We then summarize how the
theory of \cite{GuJa13} applies and state a well-posedness result for martingale solutions
in sufficiently hyperviscous cases. In \autoref{sec:details} we then collect
details and computations completing the proof of our main results.

\renewcommand{\abstractname}{Acknowledgements}
\begin{abstract}
	The authors would like to thank Franco Flandoli and Massimiliano Gubinelli 
	for many useful suggestions and relevant conversations,	
	and Martin Saal for introducing them to the literature on Primitive Equations. 
\end{abstract}

%%%%%%%%%%%%%%%%%%%%%%%%%%%%%%%%%%%%%%%%%%%%%%%%%%%%%%%%%%%%%%%%%%%%%%%%%%%%%%%%%%%%%%%%%%%%
\section{Vorticity Formulation and Conservation Laws}\label{sec:vortform}

The model under consideration in the remainder of the paper is the following stochastic
PDE in the space domain $D=[0,2\pi]^2\ni(x,z)$,
\begin{equation}\label{eq:2dspe}
\begin{cases}
\de_t v + v \de_x v + w \de_z v + \de_x p = -(-\Delta)^\theta v + \eta, \\
\de_z p = 0, \\
\de_x v + \de_z w = 0.
\end{cases}
\end{equation}
Here, $v=v(t,x,z)$ is the horizontal velocity, $w=w(t,x,z)$ is the vertical velocity,
$p=p(t,x)$ is the pressure. The parameter $\theta$ and the additive Gaussian noise $\eta$
will be specified below. The unknown fields $v,w$ are subject to the following boundary conditions:
\begin{equation}\label{eq:bc}
	\begin{cases}
	w=0, & \mbox{ if } z=0,2\pi,\\\
	v=0, & \mbox{ if } x=0,2\pi,\\
	\de_z v=0, & \mbox{ if } z=0,2\pi.
	\end{cases}
\end{equation}
The first two lines impose \emph{impermeability} of the boundary;
the third one is called a \emph{free boundary condition} for the surface and the bottom of $D$.
Before moving on, we discuss another possible choice in the next paragraph.

\subsection{On Physically Realistic Boundary Conditions}\label{ssec:physrel}

While free boundary conditions are suited to describe interfaces between fluids
such as the ocean surface, they can not be used to model a solid boundary such as the ocean bottom.
Instead, one should consider a no-slip boundary condition, leading to a different set of conditions:
\begin{equation}\label{eq:physbc}
\begin{cases}
w=0, & \mbox{ if } z=0,2\pi,\\\
v=0, & \mbox{ if } x=0,2\pi,\\
v=0, & \mbox{ if } z=0,\\
\de_z v=0, & \mbox{ if } z=2\pi.\\
\end{cases}
\end{equation}
In other words, we are assuming that the full velocity field $(v,w)$ vanishes on the bottom side.

We prefer the choice \eqref{eq:bc} since
Laplace operator can be diagonalised on functions satisfying that set of boundary conditions.
This is not true when we consider Dirichlet boundary at the bottom, since
the eigenvalue problem is overdetermined. 
In that case, Fourier analysis can still be carried through with an orthonormal basis
differing from usual trigonometric functions, see \cite[Section 6]{BrKaLe05}.
We also refer to \cite{GHTe11} for further discussion on boundary condition,
and conclude the paragraph observing that one can reduce conditions \eqref{eq:physbc} to the ones \eqref{eq:bc}. 

Assume that $(v,w)$ is a smooth solution of \eqref{eq:2dspe} on $D$ satisfying \eqref{eq:physbc},
for simplicity in the case $\eta=0$. Then, if we extend the solution to the doubled domain 
$\tilde D=[0,2\pi]\times [-2\pi,2\pi]$ so that $v,w$ are odd functions in the $z$ direction,
we have obtained a solution of \eqref{eq:2dspe} on $\tilde D$ satisfying \eqref{eq:bc}.
The size and aspect ratio of the domain is in fact irrelevant in our discussion.

\subsection{Vorticity Formulation}\label{ssec:vortform}
Let us first assume to be dealing with smooth solutions of \eqref{eq:2dspe},
driven by a smooth deterministic $\eta$.
The aim is to derive an equivalent formulation of the model
in terms of the only scalar field {vorticity} $\omega=\de_z v$,
on which we will focus the remainder of our discussion.

First of all, let us notice that $v$ must always have zero average in the $z$ direction,
since the incompressibility equation $\de_x v+\de_z w=0$ and boundary conditions imply, for all $x\in[0,2\pi]$:
\begin{equation*}
	\de_x \int_0^{2\pi} v(x,z')dz'=\int_0^{2\pi} \de_x v(x,z')dz'=-\int_0^{2\pi} \de_z w(x,z')dz'=0,
\end{equation*}
from which it follows
\begin{equation}\label{eq:vzeroaverage}
	\int_0^{2\pi} v(x,z')dz'=\int_0^{2\pi} v(0,z')dz'=0.
\end{equation}
Because of this, the solution $A(v)$ of the linear problem
\begin{equation*}
	\begin{cases}
	-\de_z^2 A(v)(x,z)=v(x,z), &(x,z)\in [0,2\pi]\times (0,2\pi),\\
	A(v)(x,z)=0, &z=0,2\pi,
	\end{cases}
\end{equation*}
is well defined for all $v$ satisfying our hypothesis. 
Another property of solutions $(v,w)$ holding independently of time is that 
$w$ is a \emph{diagnostic variable}, \emph{i.e.} it is completely determined by $v$:
\begin{equation*}
	w (x,z) = w(x,0) -\int_0^z \de_x v(x,z') dz'=\de_x\de_z A(v)(x,z).
\end{equation*}

Neglecting for a moment boundary conditions,
equations \eqref{eq:2dspe} can thus be rewritten in terms of only $v,p$ by
\begin{equation*}
\begin{cases}
\de_t v + v \de_x v +\de_x\de_z A(v) \de_z v + \de_x p =-(-\Delta)^\theta v + \eta,\\
\de_z p = 0.
\end{cases}
\end{equation*}
The system is then further simplified by considering the equation for vorticity $\omega=\de_z v$,
which does not involve the pressure $p$:
\begin{equation}\label{eq:vorticityformulation}
	\de_t \omega+ \nabla^\perp A(\omega)\cdot \nabla\omega=-(-\Delta)^\theta \omega + \de_z \eta,
\end{equation}
where $\nabla^\perp=(-\de_z,\de_x)$.
Notice that $v$ is completely determined by its partial derivative $\de_z v$ and
the zero average condition \eqref{eq:vzeroaverage}, so \eqref{eq:vorticityformulation}
is equivalent to \eqref{eq:2dspe}. Let us also remark that $A(\omega)$ is well-defined
since $\omega$ has zero average in the $z$ direction, and that $A$ --to be rigorously
defined below as an operator on function spaces-- commutes with derivatives.

Let us briefly discuss boundary conditions for $\omega$. Conditions on $v$ immediately prescribe 
$\omega(x,0)=\omega(x,2\pi)=0$ for $x\in[0,2\pi]$ and,
moreover, since $v$ is constant along the $z$ direction at $x=0,2\pi$,
we also have $\de_z v(0,z)=\de_z v(2\pi,z)=0$ for all $z\in [0,2\pi]$:
overall $\omega$ must vanish on $\de D$.
The condition $w=0$ on $z=0,2\pi$ is not as easy to translate into a condition for $\omega$,
but we will bypass the issue with our Fourier series approach below.
It is worth noticing, however, that it is because of the boundary condition on $w$
that $A(\omega)$ is well defined.

\begin{rmk}
	The relation between boundary conditions for $(v,w)$ and $\omega$ is thoroughly discussed
	in \cite{BrKaLe05} in the setting of \autoref{ssec:physrel}.
\end{rmk}

To conclude the paragraph, let us observe that thanks to the driving vector field $\nabla^\perp A(\omega)$
being Hamiltonian, smooth solutions of the hydrostatic Euler equation \eqref{eq:2dhe} in vorticity form,
\begin{equation}\label{eq:2dhevort}
	\de_t \omega+ \nabla^\perp A(\omega)\cdot \nabla\omega=0,
\end{equation}
with $(v,w)=\nabla^\perp A(\omega)$ satisfying boundary conditions, preserve the quadratic observable $\int_D \omega^2 dxdz$. Quite remarkably, this feature is peculiar to the two-dimensional case, since the quantity $\omega$ does not seem to have a counterpart in higher dimensions.

\subsection{A Rigorous Functional Analytic Setting}\label{ssec:analfun}
As we described above, we are not interested in regular solutions of \eqref{eq:2dspe},
but rather to singular, distributional regimes. It is thus convenient to
encode in Fourier series the boundary conditions, and then set up
our results in distribution spaces defined by means of Fourier expansions.

The general Fourier series expansion of a smooth function $\omega$ on $D$
such that $A(\omega)$ is well-defined and $(v,w)=\nabla^\perp A(\omega)$ satisfy boundary conditions \eqref{eq:bc} is 
\begin{equation*}
\omega(x,z)=\sum_{k\in\N^2_0} \hat \omega_k e_k(x,z),
\quad e_k(x,z)=\frac1\pi \sin(k_1 x)\sin(k_2 z),
\end{equation*}
where the $e_k$'s form an orthonormal set in $L^2(D)$, $\hat \omega_k$ are the Fourier coefficients of $\omega$
and $k=(k_1,k_2)\in \N^2_0=(\N\setminus\set{0})^2$.
We will denote
\begin{equation*}
	\S=\set{\omega=\sum_{k\in\N^2_0} \hat \omega_k e_k: \forall p\in\R\, \sum_{k \in \N^2_0} |k|^p \abs{\hat\omega_k}<\infty}.
\end{equation*}
Equivalently, $\S$ is the space of smooth functions $\omega$ on $D$ belonging to the domain of $A$
and  such that $(v,w)=\nabla^\perp A(\omega)$ satisfies the boundary conditions \eqref{eq:bc}. 
We then denote by $\S'$ its dual space,
represented by Fourier series whose coefficients grow at most polynomially.
Brackets $\brak{\cdot,\cdot}$ will denote duality couplings between functions and distributions
\begin{equation*}
	\brak{f,g}=\sum_{k\in\N^2_0} \hat f_k \hat g_k,
\end{equation*}
defined whenever the right-hand side converges. Let us also introduce, for $m\in \N$,
the projection onto the linear space of functions generated by $e_k$ with $|k|^2=k_1^2+k_2^2\leq m^2$,
\begin{equation*}
	\pi_n: \S'\to \S, \quad \omega\mapsto \pi_m\omega=\sum_{\substack{k \in \N^2_0, \\|k|\leq m}} \hat \omega_k e_k.
\end{equation*}
Following \cite{GuJa13}, we set up our analysis on the Banach spaces
\begin{equation*}
	\F L^{p,\alpha}=\set{\omega\in\S': 
		\norm{\omega}^p_{\F L^{p,\alpha}}=\sum_{k\in\N^2_0} |k|^{\alpha p}|\hat\omega_k|^p<\infty},
	\quad \alpha\in\R,p\geq 1,
\end{equation*}
and their $p=\infty$ version with $\norm{\omega}_{\F L^{\infty,\alpha}}=\sup_{k\in\N^2_0} |k|^\alpha |\hat\omega_k|$. 

Moving to the Fourier expression of the dynamics \eqref{eq:2dspe}, the crux is clearly the nonlinear term,
whose Fourier expansion is given by
\begin{equation}\label{eq:fourierb}
\nabla^\perp A(\omega)\cdot\nabla\omega=B(\omega)=\sum_{k\in\N^2_0} B_k(\omega)e_k,
\quad
B_k (\omega) =\sum_{h\in\Z^2_0} \hat\omega_h \hat\omega_{k- h} \frac{k \cdot h^{\perp}}{h_2^2},
\end{equation}
where $\Z_0^2 = (\Z\setminus\{0\})^2$ and, for $h=(h_1,h_2) \in \Z_0^2$, $\hat\omega_h = \text{sign}(h_1 h_2)\hat\omega_{(|h_1|,|h_2|)}$.

With vorticity formulation at hand, the difficulty inherent to the nonlinear term
is now apparent: looking at the $z$ component of the divergence-less vector field
$\nabla^\perp A(\omega)$, the loss of one $\de_x$ derivative is not compensated
by the gain of one $\de_z$ derivative. Indeed, such unbalance marks the difference
between \eqref{eq:2dspe} and 2-dimensional SNS, which is especially evident in the
Fourier series expansion \eqref{eq:fourierb}.

\subsection{Gaussian Invariant Measures and Driving Noise}
Referring to \cite{DPZa14}, we now introduce the stochastic analytic tools we will employ below.

Invariance of $S(\omega)=\frac12\int_D \omega(x,z)^2 dxdz$ for \eqref{eq:2dhe}
suggests that existence of an invariant \emph{Gibbs measure} formally defined by
\begin{equation}\label{eq:enstrophymeasure}
d\mu(\omega)=\frac1Z e^{- S(\omega)}d\omega.
\end{equation} 
Since $S$ is quadratic, \eqref{eq:enstrophymeasure} can be understood as a Gaussian measure on $\S'$
with covariance operator $\id$, a multiple of \emph{space white noise} on $D$.
In other words, $\mu$ is the law of the centred Gaussian process $\chi$ indexed by $\F L^{2,0}$ with covariance
\begin{equation*}
\expt{\chi(f)\chi(g)}= \brak{f,g}, \quad f,g\in \F L^{2,0}.
\end{equation*}
Such $\mu$ can be interpreted as the law of a random distribution 
supported on all $\F L^{2,\alpha}$ with $\alpha<-1$, the spaces into which the reproducing kernel
Hilbert space $\F L^{2,0}$ has Hilbert-Schmidt embedding. 
Although a fixed realisation of the random field $\chi$
is only a distribution, couplings $\brak{f,\chi}=\chi(f)$ for $f\in \F L^{2,0}$ 
are defined as random variables in $L^2(\mu)$ (It\=o integrals).

Another equivalent formulation is in terms of infinite products: formally expanding $S$ by
Parseval formula, we can write
\begin{equation*}
d\mu(\omega)=\prod_{k\in\N^2_0} \pa{\frac1{\sqrt{2\pi}}e^{-\frac12|\hat\omega_k|^2}d\hat\omega_k},
\end{equation*}
that is, under $\mu$ the Fourier coefficients $\hat \omega_k$ are independent identically distributed
standard Gaussian variables.
As a consequence, for all $\alpha<0$, $\mu$ is supported by $\F L^{\infty,\alpha}$.

Looking at the laws of Fourier components under $\mu$ it is also clear why
under this measure equations \eqref{eq:2dhe} and \eqref{eq:maineq} are \emph{singular}:
the series defining a single coefficient $B_k(\omega)$ of the vector field
diverges almost surely under $\mu$.
On the other hand, the expected value under $\mu$ of each summand in the series defining $B_k(\omega)$
vanishes, which is a formal but suggestive argument supporting the invariance of $\mu$.
In fact, the argument becomes rigorous when considering Galerkin truncations of $B$,
and we will make essential use of this in the following.
 
The \emph{space-time} analogue of $\mu$, which we will use to define the stochastic forcing
for \eqref{eq:2dspe}, can be defined in two equivalent ways.
First, we can consider the centred Gaussian field $\xi$ indexed by $L^2([0,T], \F L^{2,0})$,
with $T\in [0,\infty]$, whose covariance is given by
\begin{equation*}
	\expt{\xi(\Phi)\xi(\Phi')}=\brak{\Phi,\Phi'}_{L^2([0,T], \F L^{2,0})}.
\end{equation*}
When coupled with test functions of the form $\one_{[0,t]}(s)\phi(x,z)$, $\phi\in \S$,
$\xi$ can be regarded as the cylindrical Wiener process $W_t$ on $\F L^{2,0}$:
\begin{equation*}
\brak{\xi, \one_{[0,t]}\phi}=\brak{W_t,\phi}=\sum_{k\in\N^2_0} \hat\phi_k \beta^k_t,
\end{equation*}
the latter part being the usual Karhunen-Lo\`eve decomposition with 
$(\beta^k_t)_{k\in \N^2_0}$ independent standard Wiener processes.

For all $\theta>0,\nu>0$, the Gaussian measure $\mu$ is the unique, 
ergodic invariant measure of the infinite-dimensional Langevin's dynamics
\begin{equation}\label{eq:langevin}
	\de_t X=-\nu (-\Delta)^\theta X + \sqrt{2\nu}(-\Delta)^{\theta/2}\xi,
\end{equation}
which can be interpreted, by means of Fourier decomposition, as the system of independent one-dimensional SDEs
\begin{equation*}
	d\hat X_k=-\nu |k|^{2\theta} \hat X_k dt + \sqrt{2\nu}|k|^{\theta} d\beta^k_t, \quad k \in \N_0^2.
\end{equation*}
For the sake of simplicity, and without loss of generality, we will set $\nu=1$ in the following.
Let us conveniently introduce a symbol for the Generator of the dynamics \eqref{eq:langevin}:
first we define cylinder functionals on $\S'$ by
\begin{equation*}
	\cyl = \set{F\in L^2(\mu): \, F(\omega)=f(\hat\omega_{k_1},\dots \hat\omega_{k_r}), 
		\, f\in C^\infty(\R^r), \,k_1,\dots k_r\in \N^2_0, r\in\N},
\end{equation*} 
and for $F\in\cyl$ we denote
\begin{equation}\label{eq:ougenerator}
	\L_\theta F(\omega)= \sum_{i=1}^r |k_i|^{2\theta} \pa{-\hat\omega_{k_i} \de_i f+\de_i^2 f}.
\end{equation}
Let us also introduce the \emph{carr\'e du champ} of the diffusion operator $\L_\theta$:
for $F,G\in\cyl$,
\begin{equation}\label{eq:carreduchamp}
	\E_\theta(F,G)(\omega)=\sum_{i=1}^n |k_i|^{2\theta} \de_i f\de_i g,
\end{equation}
which satisfies the Gaussian integration by parts formula
\begin{equation*}
	\expt[\mu]{F \L_\theta G}=-\expt[\mu]{\E_\theta(F,G)}.
\end{equation*}

The above arguments finally lead us to consider the combination of dynamics \eqref{eq:2dhe} and \eqref{eq:langevin}
as a SPDE preserving $\mu$:
\begin{equation}\label{eq:maineq}
	\de_t \omega+ \nabla^\perp A(\omega)\cdot \nabla\omega
	=-(-\Delta)^\theta  \omega + \sqrt{2}(-\Delta)^{\theta/2}\xi.
\end{equation}
As already noticed, the nonlinear part of the dynamics
is not well defined for functions $\omega$ in the regularity regime dictated by $\mu$,
or rather, it can be given a rigorous meaning only by exploiting cancellations
due to the structure of the stochastic equation as a whole. 

\begin{rmk}
	In terms of $v$, the latter equation reads
	\begin{equation*}
	\begin{cases}
	\de_t v + v \de_x v +\de_x\de_y A(v) \de_y v + \de_x p 
	=- (-\Delta)^\theta v +\de_z \sqrt{2} (-\Delta)^{\theta/2}\xi, \\
	\de_z p = 0.
	\end{cases}
	\end{equation*}
	The forcing term should have white noise regularity in $x$, and Brownian regularity in $y$,
	although the covariance structure is a nontrivial copula of the two. 
\end{rmk}

\begin{rmk}
	Just as in the case of 2D Euler or stochastic Navier-Stokes equations,
	the invariant measure associated to enstrophy is not able to describe
	peculiar features of the fluid-dynamic model, such as turbulence phenomena.
	In fact, such measures are preserved
	by any flow of measure-preserving diffeomorphisms of the domain,
	among which the Euler flow is a very distinguished case.
	Energy ensembles should be in fact more relevant, but they are supported
	on quite larger distribution spaces.
\end{rmk}

%%%%%%%%%%%%%%%%%%%%%%%%%%%%%%%%%%%%%%%%%%%%%%%%%%%%%%%%%%%%%%%%%%%%%%%%%%%%%%%%%%%%%%%%%%%%%
\section{Regularisation by Noise in Hyperviscous Regimes} 

In this section we outline how the solution theory of \cite{GuJa13}
(known as \emph{Energy Solutions} theory in the context of stochastic Burgers and KPZ equations)
applies to our model in a sufficiently hyperviscous regime.
Computations differ from that work only by small details: we collect them in 
the last section for the sake of completeness, and in the present one we only
recall the core ideas.

\subsection{Controlled Processes and Martingale Solutions}
We recall the notion of controlled process from \cite{GuJa13}.

\begin{definition}\label{def:controlledprocess}
	For $\theta\geq 0$ and $T> 0$ we define the space $\RR_{\theta,T}$ of stochastic processes
	with trajectories of class $C([0,T],\S')$ such that any $\omega\in \RR_{\theta,T}$ satisfies:
	\begin{enumerate}
		\item $\omega$ is stationary and for any $t\in [0,T]$, $\omega_t\sim \eta$;
		\item there exists a stochastic process $\A$ with trajectories $C([0,T],\S')$
		starting from $\A_0=0$ and with null quadratic variation such that,
		for any $\phi\in\S$,
		\begin{equation*}
			\brak{\phi,\omega_t}-\brak{\phi,\omega_0}
			+\int_0^t \brak{(-\Delta)^\theta\phi,\omega_s}ds-\brak{\phi,\A_t}=M_t(\phi)
		\end{equation*}
		is a martingale with respect to the filtration of $\omega$, and it has quadratic variation
		$\bra{M(\phi)}_t=2t \norm{(-\Delta)^{\theta/2}\phi}^2_{\F L^{2,0}}$;
		\item the reversed process $\tilde \omega_t=\omega_{T-t}$ satisfies condition (2)
		with $\tilde \A_t=-\A_{T-t}$.
	\end{enumerate} 
\end{definition}

Notice that in fact elements of $\RR_{\theta,T}$ are the couples $(\omega,\A)$.
The forward and backward martingale equations defining the class $\RR_{\theta,T}$ 
allow to obtain good \emph{a priori} estimates for nonlinear functionals of controlled process,
in a procedure by now commonly known as \emph{It\=o trick}, especially in literature
related to regularisation by noise techniques, see \cite{Fl10,FlGuPr10,BeFlGuMa19}.

In the next paragraph we detail how the It\=o trick produces good estimates on Galerkin approximations of \eqref{eq:maineq}:
the idea behind \autoref{def:controlledprocess} is to collect the features of those approximants
allowing such estimates, to form a class of processes on which the nonlinear term of \eqref{eq:maineq} is defined.
In \autoref{sec:details} we will prove the following:

\begin{lemma} \label{lem:nonlinearity}
	Let $\theta > 2$, $T> 0$ and $\omega\in \RR_{\theta,T}$. Then for every $\zeta < -1$
	\begin{equation*}
\lim_{m\to \infty} \int_0^t B(\pi_m \omega_s)ds
\end{equation*}
exists as a limit in $C([0,T],\F L^{\infty,\zeta})$.
	We denote by $\int_0^t B(\omega_s)ds$ the limiting process.
\end{lemma} 

The latter lemma shows that the nonlinear functional $B(\omega)$ can be defined for $\omega\in\RR_{\theta,T}$
as a distribution in both space \emph{and} time.
Let us observe that Fourier truncation $\pi_m$ in \autoref{lem:nonlinearity} can in fact be replaced
with a large class of mollifiers, the limit being independent of such choice:
for the sake of keeping the exposition simple, we refrain from going into details. 

We can now give a notion of martingale solution to \eqref{eq:maineq}.

\begin{definition}\label{def:martsolution}
	Let $\theta> 2$, $T> 0$ and $\omega\in \RR_{\theta,T}$. We say that $\omega$
	is a \emph{martingale solution} to \eqref{eq:maineq} if it holds almost surely, for any $t\in[0,T]$,
	\begin{equation*}
		\A_t=\int_0^t B(\omega_s)ds.
	\end{equation*}
	The solution is \emph{pathwise unique} if, for any two controlled processes 
	$\omega,\tilde\omega\in\RR_{\theta,T}$ defined on the same probability space,
	satisfying conditions (2) and (3) of \autoref{def:martsolution} with the same martingales and
	with $\omega_0=\tilde\omega_0$ almost surely, then almost surely, for all $t\in [0,T]$, $\omega_t=\tilde\omega_t$.
\end{definition}

The following is the main result of the paper: its proof will be given in \autoref{sec:details}.

\begin{thm} \label{thm:main}
	Let $T>0$. For any $\theta>2$ there exists a solution to \eqref{eq:maineq} in the sense of \autoref{def:martsolution}.
	Moreover, for $\theta>3$ the solution is pathwise unique. 
\end{thm}

\subsection{Galerkin Approximation and the It\=o Trick}
Let us introduce approximating processes $(\omega^m)_{m\in\N}$ by their Fourier coefficients dynamics:
for $k\in\N^2_0$,
\begin{align}\label{eq:galerkin}
	d \hat\omega_k^m &= B^m_k(\omega^m) dt - |k|^{2\theta} \hat\omega^m_k dt + \sqrt{2}|k|^\theta d\beta^k_t,
\end{align}
where $B^m(\omega)=\pi_m B(\pi_m\omega)$, and $\omega^m_0\sim \mu$. 
The vector field $B^m$ satisfies
\begin{align}\label{eq:divergencefree}
\div_\mu B^m(\omega)
&= \div_\mu \sum_{\substack{k \in \N^2_0, \\|k|\leq m}}\sum_{\substack{h \in \Z^2_0, \\|h|\leq m}} 
\hat\omega_h \hat\omega_{k- h} \frac{k \cdot h^{\perp}}{h_2^2} e_k\\ \nonumber
&= \sum_{\substack{k \in \N^2_0, \\|k|\leq m}}\sum_{\substack{h \in \Z^2_0, \\|h|\leq m}} \pa{\de_{\hat\omega_k}\pa{\hat\omega_h \hat\omega_{k- h}}
-\hat\omega_h \hat\omega_{k- h} \hat\omega_k} \frac{k \cdot h^{\perp}}{h_2^2}=0.
\end{align}
As a consequence, \eqref{eq:galerkin} has a unique, (probabilistically) strong, global in time solution since 
$\mu$ is preserved by the linear part of the dynamics, and thus \cite[Theorem 3.2]{Cr83} applies.
In the following, we denote by $\PP^m_\mu$ the law of $\omega^m$ in $C(\R_+,\S')$.

By It\=o formula, for any cylinder function $F\in\cyl$,
$F(\omega) = f(\hat\omega_{k_1},\dots, \hat \omega_{k_n})$, it holds
\begin{equation*}
	d F(\omega^m)=\L_\theta F(\omega^m)dt+ \G^m F(\omega^m)dt
	+ \sum_{i=1}^n \de_i f(\hat\omega_{k_1}^m,\dots,\hat \omega_{k_n}^m) \sqrt{2}|k_i|^\theta d\beta^{k_i}_t,
\end{equation*}
where $\L_\theta$ is defined in \eqref{eq:langevin} and
\begin{equation*}
	\G^m F(\omega)=\sum_{i=1}^n \de_i f(\hat\omega_{k_1},\dots, \hat\omega_{k_n}) B^m_{k_i}(\omega) dt.
\end{equation*}
In other words, the process
\begin{equation}\label{eq:mart1}
	M^{F,m}_t= F(\omega_t^m) - F(\omega_0^m) - \int_0^t \L_\theta F (\omega^m_s) ds
 		- \int_0^t \G^m F (\omega^m_s) ds
\end{equation}
is a martingale with quadratic variation
\begin{equation*}
[M^{F,m}]_t = 2 \int_0^t \sum_{i=1}^n |k_i|^{2\theta} \pa{\de_i f(\hat\omega_{k_1}^m,\dots,\hat \omega_{k_n}^m)}^2 ds 
= 2 \int_0^t \E_\theta(F)(\omega_s^m) ds.
\end{equation*}

Let us point out that, thanks to the hydrodynamic form of the nonlinearity, 
$\G^m$ is a skew-symmetric operator with respect to $\mu$:
indeed, since
\begin{equation*}
\brak{\omega,B^m(\omega)}=\brak{\omega, \pi_m(\nabla^\perp A(\pi_m \omega)\cdot \nabla\pi_m\omega)}=0,
\end{equation*}
Gaussian integration shows that
\begin{equation*}
\expt[\mu]{F\G^m G}=-\expt[\mu]{G\G^m F}, \quad \forall F,G\in\cyl.
\end{equation*}

Let us then consider the reversed process $\tilde \omega^m_t=\omega^m_{T-t}$, for a fixed time horizon $T>0$:
$\tilde \omega^m$ is a Markov process whose generator is the adjoint of the one of $\omega^m$,
that is $\L_\theta-\G^m$. The process
\begin{equation}\label{eq:mart2}
	\tilde M^{F,m}_t= F(\tilde \omega_t^m) - F(\tilde \omega_0^m) - \int_0^t \L_\theta F (\tilde \omega^m_s) ds
	- \int_0^t \G^m F (\tilde \omega^m_s) ds
\end{equation}
is thus another martingale with quadratic variation $2 \int_0^t \E_\theta(F)(\omega_s^m) ds$.
To sum up, we have shown that $\omega^m$ is a controlled process in the sense of \autoref{def:controlledprocess}.

The trick is now to sum the martingale identities \eqref{eq:mart1}, \eqref{eq:mart2} 
for $\omega^m$ and $\tilde{\omega}^m$: in doing so the nonlinear skew symmetric part, 
together with boundary terms, is canceled,
leaving us with martingales term and the symmetric Ornstein-Uhlenbeck generator,
\begin{align*}
	\tilde{M}^{F,m}_{T-t} - \tilde{M}^{F,m}_T - M^{F,m}_t = 2 \int_0^t  \L_\theta F (\omega^m_s) ds.
\end{align*}
Burkholder-Davis-Gundy inequality thus provides, together with stationarity of $\omega^m$,
the following powerful estimate: for $p\geq 1$ there exists a constant $C_p>0$ only depending on $p$ such that
for all $F\in\cyl$
\begin{equation} \label{eq:itotrick}
	\expt[\PP^m_\mu]{\sup_{t \in [0,T]} \abs{\int_0^t \L_\theta F (\omega_s^m) ds}^p}
		\leq C_p \sqrt T \expt[\mu]{\abs{\E_\theta F}^{p/2}}.
\end{equation}

\begin{rmk}
	As already observed, \autoref{def:controlledprocess} actually collects the elements we used to
	establish \eqref{eq:itotrick}; indeed, the latter holds more generally for any controlled process $\omega\in\RR_{\theta,T}$.
\end{rmk}

Inequality \eqref{eq:itotrick} provides good estimates on time integrals of observables for $\omega^m$,
provided that we are able to solve a Poisson equation in Gaussian space.
The main aim are clearly bounds to establish the limit in \autoref{lem:nonlinearity},
which can be obtained by means of \eqref{eq:itotrick} by solving
\begin{equation*}
	\L_\theta H^m_k(\omega)=B^m_k(\omega).
\end{equation*}
Since $B^m_k(\omega)$ belongs to the second chaos in the Wiener chaos decomposition of $L^2(\mu)$,
and since $\L_\theta$ is diagonalised by such decomposition, it is easy to obtain the explicit solution
\begin{equation}\label{eq:soluzh}
	H_k^m (\omega) = - \chi_{\{|k| \leq m\}} \sum_{\substack{h,\ell \in \Z^2_0, \\ h + \ell = k \\ |h|,|\ell|\leq m}}
		\hat\omega_h \hat\omega_{\ell} \frac{\ell \cdot h^{\perp}}{h_2^2 (|h|^2 +|\ell|^2)^{\theta}}.
\end{equation}
The computation is completely analogous to \cite[Section 3]{GuJa13}, to which we refer.
With the latter expression at hand, one only needs to estimate moments of $\E_\theta(H_k^m)$:
we report such computation in the next section, together with some variants from which
\autoref{thm:main} follows.

%%%%%%%%%%%%%%%%%%%%%%%%%%%%%%%%%%%%%%%%%%%%%%%%%%%%%%%%%%%%%%%%%%%%%%%%%%%%%%%%%%%%%%%%%%%%%
\section{Proof of Main Result} \label{sec:details}

We complete in this Section the proof of \autoref{thm:main}.
First, by means of the It\=o trick estimate \eqref{eq:itotrick}
we obtain bounds on Galerkin approximations: the last two paragraphs are then
devoted to existence and uniqueness of martingale solutions. 
In this section, the symbol $\lesssim$ denotes inequality up to a positive multiplicative
constant uniform in the involved parameters.
 
\subsection{Controlling the Nonlinear Term} \label{subsec:apriori}

We start from the expression \eqref{eq:soluzh} for $H^m_k$
to obtain estimates on the nonlinear term in \eqref{eq:maineq}. 
By definition of $\E_\theta$, \eqref{eq:carreduchamp}, one has
\begin{equation*}
\E_\theta(H^m_k)(\omega)=\chi_{\{|k| \leq m\}} \sum_{\substack{h\in \Z^2_0,\\|h|\leq m}} |h|^{2\theta} \abs{\frac{2(k-h) \cdot h^{\perp}}{h_2^2(|k-h|^2+|h|^2)^\theta} \hat\omega_{k-h}}^2,
\end{equation*}
therefore, taking expectation with respect to $\mu$, for $|k|\leq m$
\begin{align*}
\mathbb{E}_\mu \left[\E_\theta(H^m_k) \right] 
&= \sum_{\substack{h\in \Z^2_0,\\|h|\leq m}}  |h|^{2\theta} \abs{\frac{2(k-h) \cdot h^{\perp}}{h_2^2(|k-h|^2+|h|^2)^\theta}}^2 \\
&\lesssim \sum_{\substack{h\in \Z^2_0,\\|h|\leq m}}   \frac{|k|^2 |h|^{2+2\theta}}{|k-h|^{4\theta}+|h|^{4\theta}}
\lesssim  \sum_{\substack{h\in \Z^2_0,\\|h|\leq m}}   \frac{|k|^2 |h|^{2}}{|k-h|^{2\theta}+|h|^{2\theta}}.
\end{align*}
Now we use the fact that, for $\theta > 2$,
\begin{equation*}
\sum_{h \in \Z^2_0} \frac{| h |^2}{| k - h |^{2 \theta}+| h |^{2 \theta}} \lesssim
| k |^{4 - 2 \theta}
\end{equation*}
(see \cite[Lemma 16]{GuTu19}) to deduce the following estimate uniformly in $m$:
\begin{equation*}
\mathbb{E}_\mu \left[\E_\theta(H^m_k) \right]  \lesssim | k |^{6 - 2 \theta}.
\end{equation*}
Similarly, increments are controlled by
\begin{equation*}
\sup_{n>m} \mathbb{E}_\mu [ \mathcal{E}_{\theta} (H_k^n-H_k^m)]
	\lesssim | k |^2 m^{ 4-2\theta} .
\end{equation*}
With these estimates at hand, by means of \eqref{eq:itotrick} and Gaussian hypercontractivity,
one can prove the following estimates on the nonlinear term of \eqref{eq:maineq}.

\begin{lemma} \label{lem:apriori1}
	Let $G^m_t \coloneqq \int_0^t B(\pi_m \omega_s)ds$ and
	${\mathbb{P}^m_\mu}$ be the distribution of the stationary solution of \eqref{eq:galerkin} described above.  
	For any $n>m$ we have the following estimates:
	\begin{gather}
	\left\| \sup_{t \in [0,T]} \left( G^{m}_t\right)_k \right\|_{L^p({\mathbb{P}^m_\mu})} \lesssim |k|^{3 - \theta} T^{1/2}, \label{eq:apriori2}\\
	\left\| \sup_{t \in [0,T]} \left( G^{n}_t\right)_k - \left( G^{m}_t\right)_k \right\|_{L^p({\mathbb{P}^m_\mu})} \lesssim |k| T^{1/2} m^{2-\theta}. \label{eq:apriori3}
	\end{gather}
\end{lemma}

\begin{lemma} \label{lem:apriori2}
	Let $\tilde{G}^m_t \coloneqq \int_0^t e^{-(t-s)(-\Delta)^\theta} B(\pi_m \omega_s)ds$ and 
	${\mathbb{P}^m_\mu}$ as above. 
	For any $m$ fixed, $n>m$, $s,t \in [0,T]$, $s<t$,  we have the following estimates:
	\begin{gather}
	\left\| \sup_{t \in [0,T]} \left( \tilde{G}^{m}_t\right)_k \right\|_{L^p({\mathbb{P}^m_\mu})} \lesssim |k|^{3 - 2\theta}, \label{eq:apriori5}\\
	\left\| \sup_{t \in [0,T]} \left( \tilde{G}^{n}_t\right)_k - \left( \tilde{G}^{m}_t\right)_k \right\|_{L^p({\mathbb{P}^m_\mu})} \lesssim |k|^{-1} m^{4-2\theta}, \label{eq:apriori6}\\
	\sup_{m}\left\| \left( \tilde{G}^{m}_t\right)_k - \left( \tilde{G}^{m}_s\right)_k \right\|_{L^p({\mathbb{P}^m_\mu})} \lesssim |k|^{3-2\theta+2\varepsilon\theta} (t-s)^{\varepsilon},\label{eq:aprioricont}
	\end{gather}
	where the last inequality is meant to hold for $\varepsilon>0$ small enough.
\end{lemma}

Proofs of the previous estimates follow along the lines of Lemma 5, 
Lemma 6 and Corollary 1 of \cite{GuJa13}, so we refrain from repeating them here.

\subsection{Existence for $\theta>2$}
We first prove \autoref{lem:nonlinearity}, which gives a meaning to the nonlinear term of \eqref{eq:maineq}. 
The result easily follows from \autoref{lem:apriori1}. 

\begin{proof}[Proof of \autoref{lem:nonlinearity}]
	Let $G^m_t \coloneqq \int_0^t B(\pi_m \omega_s)ds$.
	It is clear that $G^m$ is a random process with values in $C([0,T],\F L^{\infty,\zeta})$ for every $m$ and $\zeta \in \R$. 
	Since $\theta>2$, \eqref{eq:apriori3} gives for any $p$ and $n>m$
	\begin{equation*}
	\mathbb{E}_{\mathbb{P}^m_\mu} \left[ \sum_k |k|^{\zeta p} \left| \sup_{t \in [0,T]} \left( G^{n}_t\right)_k - \left( G^{m}_t\right)_k \right|^p \right] \to 0
	\end{equation*}
	as $m \to \infty$ whenever $\zeta <-2/p -1$. 
	Taking $p$ sufficiently large, for any $\zeta < -1$ we obtain the almost sure uniform 
	convergence of $G^m$ in the space $C([0,T],\F L^{\infty,\zeta})$.
\end{proof}

We are now ready to prove the first part of \autoref{thm:main}. 
The proof relies on \autoref{subsec:apriori} and Skorokhod Theorem.
 
\begin{proof}[Proof of \autoref{thm:main}, existence]
	Let us consider the mild formulation of \eqref{eq:galerkin}: 
	\begin{align} \label{eq:Galerkin}
	\omega^m_t &= e^{-t(-\Delta)^\theta} \omega_0 + \int_0^t e^{-(t-s)(-\Delta)^\theta} B^m(\omega^m_s) ds \\ \nonumber
	&\quad + \sqrt{2}(-\Delta)^{\theta/2} \int_0^t e^{-(t-s)(-\Delta)^\theta} d\beta_s,
	\end{align}
	where $\omega_0 \sim \mu$ and $\beta$ is a cylindrical Wiener process on $\F L^{2,0}$. Define 
	\begin{equation*}
	\A^m_t \coloneqq \int_0^t B^m(\omega^m_s) ds , \quad \tilde{\A}^m_t \coloneqq \int_0^t e^{-(t-s)(-\Delta)^\theta} B^m(\omega^m_s) ds. 
	\end{equation*}
	We prove that, for  $\varepsilon>0$ sufficiently small and $\zeta < -1$, the laws of the processes $\left( \omega^m , \A^m, \tilde{\A}^m, \beta \right)_m$ are tight in $C([0,T],\mathcal{X})$, where 
	\begin{equation*}
	\mathcal{X} \coloneqq \F L^{\infty,\zeta} \times \F L^{\infty,\theta - 3 - \varepsilon} \times \F L^{\infty,2\theta - 3 -\varepsilon} \times\F L^{\infty,-\varepsilon}.
	\end{equation*}
	By Borel-Cantelli theorem applied to Fourier expansions, the law $\mu$ is concentrated on $\F L^{\infty,-\varepsilon}$, 
	and the stochastic convolution takes values in $C([0,T],\F L^{\infty,\theta -\varepsilon})$ for every $\varepsilon>0$. 
	Tightness in this space is given by Fernique Theorem. 
	Tightness of $(\tilde{\A}^m)_m$ descends from \autoref{eq:aprioricont}
	and tightness of $({\A}^m)_m$ descends from \autoref{eq:apriori2}. 
	Hence, by a standard application of Prokhorov Theorem and Skorokhod Theorem, 
	we deduce the a.s. convergence, up to a subsequence and a change of the underlying probability space, 
	of $\left( \omega^m , \A^m, \tilde{\A}^m, \beta \right)$ towards some random variable 
	$\left( \omega , \A, \tilde{\A}, \beta \right)$ in $C([0,T],\mathcal{X})$ which satifies
	\begin{align*}
	\omega_t &= e^{-t(-\Delta)^\theta} \omega_0 + \tilde{\A}_t + \sqrt{2}(-\Delta)^{\theta/2} \int_0^t e^{-(t-s)(-\Delta)^\theta} d\beta_s \\
	&= \omega_0 + \int_0^t (-\Delta)^\theta \omega_s ds + \A_t + \sqrt{2}(-\Delta)^{\theta/2} \beta_t. 
	\end{align*}
	Now it is easy to check that $\omega \in \mathcal{R}_{\theta,T}$, see \cite{GuJa13} for details.
\end{proof}

\subsection{Uniqueness for $\theta>3$}

\begin{proof}[Proof of \autoref{thm:main}, uniqueness]
	We have constructed a sequence of $\omega^m$ converging a.s. to a solution $\omega$ as random variables in $C([0,T], \F L^{\infty,\zeta})$ for every $\zeta < -1$. Here we prove uniqueness, which comes from an estimate on the quantity $\pi_m(\omega^m - \omega)$ in a suitable space, where $\omega \in \mathcal{R}_{\theta,T}$ is a controlled solution to \eqref{eq:maineq} and $\omega^m$ is its Galerkin approximation defined by \eqref{eq:Galerkin}. In particular, we prove that $\pi_m(\omega^m - \omega)$ converges a.s. to zero in the space $C([0,T], \F L^{\infty,\xi})$, for suitable $\xi > \zeta$. This would conclude the proof by uniqueness at the level of the Galerkin truncations.
	
	It is easy to see that
	\begin{align*}
	\pi_m(\omega^m_t - \omega_t) \coloneqq \delta_t^m &= \int_0^t e^{-(t-s)(-\Delta)^{\theta}} (B^m(\omega^m_s) - \pi_m B(\omega_s)) ds \\
	&= \int_0^t e^{-(t-s)(-\Delta)^{\theta}} (B^m(\omega^m_s) - B^m(\omega_s)) ds \\
	&+ \int_0^t e^{-(t-s)(-\Delta)^{\theta}} (B^m(\omega_s) - \pi_m B(\omega_s)) ds \\
	&= \alpha^m_t + \gamma^m_t.
	\end{align*}
	for every $m$, and thus for every $\xi$
	\begin{align*}
	\sup_k \sup_{t \in [0,T]} |k|^\xi |(\delta^m_t)_k| \leq \sup_k \sup_{t \in [0,T]} |k|^\xi |(\alpha^m_t)_k| + \sup_k \sup_{t \in [0,T]} |k|^\xi |(\gamma^m_t)_k|.
	\end{align*}
	By \eqref{eq:apriori5}, \eqref{eq:apriori6} and interpolation, $\gamma^m$ satisfies
	\begin{equation*}
	\left\| \sup_{t \in [0,T]} \left| (\gamma^m_t)_k \right|\right\|_{L^{p}({\mathbb{P}^m_\mu})}  \lesssim |k|^{3-2\theta+\varepsilon} m^{-\varepsilon},
	\end{equation*} 
	and therefore for every $\xi< 2\theta - 3$ we have
	\begin{equation*}
	\sup_k \sup_{t \in [0,T]} |k|^\xi |(\gamma^m_t)_k| \to 0 \mbox{ a.s. for }  m \to \infty.
	\end{equation*}
	On the other hand, since 
\begin{equation*}
\left| B^m(\omega^m_s) - B^m(\omega_s) \right| \lesssim \sum_{\substack{h\in \Z^2_0,\\|h|\leq m}}  |k||h||(\omega^m_s + \omega_s)_h| |(\omega^m_s - \omega_s)_{k-h}|,
\end{equation*}	
we obtain the following bound on $\alpha^m$:
	\begin{align*}
	\sup_{t \in [0,T]} |k|^\xi |(\alpha^m_t)_k| \lesssim &|k|^\xi \sup_h \sup_{t \in [0,T]} |h|^\xi |(\delta^m_t)_h| \\
	&\times \sup_{t \in [0,T]} \left| \int_0^t e^{-(t-s)|k|^{2\theta}} \sum_{\substack{h\in \Z^2_0,\\|h|\leq m}}  |h|^{1-\xi}|k||(\omega^m_s + \omega_s)_{k-h}| ds\right|.
	\end{align*}
	If $\xi>3$ the series $\sum_{h \in \Z^2_0}|h|^{1-\xi}$ converges, therefore by H\"older inequality
	\begin{gather*}
	\left| \int_0^t e^{-(t-s)|k|^{2\theta}} \sum_{|h| \leq m} |h|^{1-\xi}|(\omega^m_s + \omega_s)_{k-h}| ds\right| \\
	\lesssim \left| \int_0^t\sum_{|h| \leq m} |h|^{1-\xi} e^{-p'(t-s)|k|^{2\theta}}  ds\right|^{1/p'}
	\left| \int_0^t \sum_{|h| \leq m} |h|^{1-\xi}|(\omega^m_s + \omega_s)_{k-h}|^p ds\right|^{1/p} \\
	\lesssim  |k|^{-2\theta/p'} \left| \int_0^t \sum_{|h| \leq m} |h|^{1-\xi}|(\omega^m_s + \omega_s)_{k-h}|^p ds\right|^{1/p}.
	\end{gather*}
	Taking $p'\to 1$ such that $\xi + 1 -2\theta/p' < 0$ and using the fact that $\omega^m$ and $\omega$ have marginals $\sim \mu$, we finally get
	\begin{equation*}
	\sup_k \sup_{t \in [0,T]} |k|^\xi |(\delta^m_t)_k| \to 0 \mbox{ a.s. for }  m \to \infty,
	\end{equation*}
	for every $3<\xi<2\theta - 3$, which corresponds to the additional contraint $\theta>3$. The proof is complete.
\end{proof}

\end{document}